\documentclass[DIV=calc, paper=letter, fontsize=11pt]{scrartcl}	 

\usepackage[]{units}
\usepackage{lipsum} 
\usepackage[english]{babel} 
\usepackage[protrusion=true,expansion=true]{microtype} 
\usepackage{amsmath,amssymb,amsfonts,amsthm} 

\usepackage[svgnames, table]{xcolor} 
\usepackage[hang, small,labelfont=bf,up,textfont=it,up]{caption} 
\usepackage{booktabs} 
\usepackage{fix-cm}	 
\usepackage[]{units}

\usepackage{subfig}

\usepackage{tabularx} 

\usepackage[colorlinks=true, linkcolor=blue, citecolor=blue, 
filecolor=blue, runcolor=blue, urlcolor = black]{hyperref}

\usepackage{fullpage}
\usepackage{sectsty} 

\usepackage{fancyhdr} 
\usepackage{lastpage} 

\usepackage{nicefrac}

\usepackage{tikz}
\usetikzlibrary{calc,patterns,decorations.pathmorphing,
decorations.markings}
\usepackage{pgfplots}

\usepackage{multicol}
\usepackage{float} 

\usepackage[utf8]{inputenc}
\usepackage[english]{babel}

\usepackage{mathtools} 

\usepackage[
backend=bibtex,
style=numeric-comp,
sorting=ynt,
maxbibnames=99,
]{biblatex}

\makeatletter
\newcommand{\multicolinterrupt}[1]{
\setcounter{tempcolnum}{\col@number}
\end{multicols}
#1%
\begin{multicols}{\value{tempcolnum}}
}
\makeatother

\newcommand{\revI}[1]{#1}

\usepackage{graphicx}
\DeclareGraphicsExtensions{.png}

\usepackage{lettrine} 
\newcommand{\initial}[1]{ 
\lettrine[lines=3,lhang=0.3,nindent=0em]{
\color{DarkGoldenrod}
{\textsf{#1}}}{}}

\usepackage{caption}

\DeclareOldFontCommand{\bf}{\normalfont\bfseries}{\mathbf}

\usepackage{graphicx}

\usepackage{titling} 

\newcommand{\HorRule}{\color{DarkGoldenrod} \rule{\linewidth}{1pt}} 

\pretitle{\vspace{-80pt} \begin{flushleft} \HorRule 
\fontsize{14}{14} \color{DarkRed} 
\selectfont} 

\title{First-Passage Time Statistics on Surfaces of General Shape: 
\\ Surface PDE Solvers using Generalized Moving Least Squares (GMLS)}

\posttitle{\end{flushleft}}

\preauthor{\vspace{-5pt} \begin{flushleft}\fontsize{12}{12} \large  
\color{DarkRed}} 
\author{B. J. Gross $^{1,2}$, P. Kuberry $^{3}$, 
P. J. Atzberger$^{1,2}$} 
\postauthor{\fontsize{12}{12} \small 
\color{Black} 
\vskip 0.5em
$[1]$ Department of Mathematics, University of California 
Santa Barbara (UCSB);  \\ 
$[2]$ Department of Mechanical Engineering, University of California 
Santa Barbara (UCSB);
\\ $[3]$ Sandia National Laboratories, Albuquerque, NM $^{**}$ 
\vskip 0.5em
atzberg@gmail.com; \hspace{0.0cm} \url{http://atzberger.org/}  
\vskip 0.5em 
\lfoot{\footnotesize * Work supported by DOE Grant ASCR PhILMS
DE-SC0019246 and NSF Grant DMS-1616353. \\ **{\scriptsize Sandia National
Laboratories is a multi-mission laboratory managed and operated by National
Technology and Engineering Solutions of Sandia, LLC., a wholly owned
subsidiary of Honeywell International, Inc., for the U.S. Department of
Energy’s National Nuclear Security Administration under contract
DE-NA0003525. This paper describes objective technical results and analysis.
Any subjective views or opinions that might be expressed in
the paper do not necessarily represent the views of the U.S. 
Department of Energy or the United States Government.} }
\end{flushleft} \vspace{-18pt} \HorRule} 
\date{}  

\newcommand{\gradX}{\nabla_{\mathbf{X}}}

\newcommand{\mb}[1]{\mathbf{#1}}
\newcommand{\bs}[1]{\boldsymbol{#1}}
\newcommand{\bsy}[1]{\boldsymbol{#1}}

\newcommand{\h}{h}

\DeclareGraphicsExtensions{.png}

\begin{document}

\maketitle 

\thispagestyle{fancy} 


\initial{W}\textbf{e develop numerical methods for computing statistics of
stochastic processes on surfaces of general shape with drift-diffusion dynamics
$d\mb{X}_t = a(\mb{X}_t)dt + \mb{b}(\mb{X}_t)d\mb{W}_t$.
\revI{We formulate descriptions of Brownian motion and general
drift-diffusion processes on surfaces.
We consider 
statistics of the form $u(\mb{x}) =
\mathbb{E}^{\mb{x}}\left[\int_0^\tau g(\mb{X}_t)dt \right] +
\mathbb{E}^{\mb{x}}\left[f(\mb{X}_\tau)\right]$ for a domain $\Omega$ and 
the exit stopping time $\tau =
\inf_t \{t > 0 \; |\; \mb{X}_t \not\in \Omega\}$, where $f,g$ are 
general smooth functions. For 
computing these statistics, we develop high-order Generalized 
Moving Least Squares (GMLS) solvers for associated 
surface PDE boundary-value problems based on 
Backward-Kolmogorov equations.  We
focus particularly on the mean First Passage Times (FPTs) given by the
case $f = 0,\, g = 1$ where $u(\mb{x}) = \mathbb{E}^{\mb{x}}\left[\tau\right]$.
We perform studies for a variety of shapes showing our methods converge with
high-order accuracy both in capturing the geometry and the surface PDE
solutions.  We then perform studies showing how statistics are influenced by the
surface geometry, drift dynamics, and spatially dependent diffusivities.}}

\setlength{\parindent}{5ex}

\section*{Introduction}
\revI{
Path-related statistics of stochastic processes, such as the 
mean First Passage Times (FPTs)~\cite{Klein1952,Redner2001,Oksendal2000},
arise in many fields, including in biology~\cite{Chou2014,
Polizzi2016, AtzbergerFPKMC2014, Ghusinga2017,Taillefumier2013,Atzberger_Tran_2021},
physics~\cite{Lindsay2017,Kusters2014, McQuarrie1990,Kells2019,Mueller1997,
Hohenegger2017}, engineering~\cite{Debnath2018,Byl2014, Byl2006},
finance~\cite{Chicheportiche2014, Bachelier1900, Hull2009, Gutlierrez2013}, and
machine learning~\cite{Callut2007, Xu2019}.  Many problems involve 
stochastic processes within manifolds where significant roles 
are played by geometric and topological contributions~\cite{Benichou2015,
Lindsay2017,Kusters2014,Ben-Naim2014,Grebenkov2018,Atzberger_Tran_2021}.  
We consider Ito processes
$\mb{X}_t$ with the drift-diffusion dynamics $d\mb{X}_t = a(\mb{X}_t)dt +
\mb{b}(\mb{X}_t)d\mb{W}_t$,~\cite{Oksendal2000,Gardiner1985}.
We formulate descriptions of Brownian motion and general
drift-diffusion processes on surfaces.} \revI{While in
principle statistics can be estimated by using stochastic numerical methods to
sample trajectories $\{\mb{X}_t\}_{0 \leq t \leq T}$, it can be
computationally expensive to reduce 
statistical sampling errors sufficiently.  Further challenges also arise 
when the process has multiple dynamical times-scales resulting in stiffness, 
or when trying to estimate statistics at many locations on the
surface~\cite{Platen1992,Gardiner1985, Mueller1997,AtzbergerFPKMC2014}.}

\revI{
In practice, for making predictions of observations and measurements, it is often 
enough to consider the class of statistics of the form
$u(\mb{x}) = \mathbb{E}^{\mb{x}}\left[\int_0^\tau g(\mb{X}_t)dt \right] +
\mathbb{E}^{\mb{x}}\left[f(\mb{X}_\tau)\right]$.  For the stochastic process,
the $\mb{X}_0 = \mb{x} \in
\Omega$ within an open domain $\Omega$ and the exit stopping time is $\tau =
\inf_t \{t > 0 \; |\; \mb{X}_t \not\in \Omega\}$.  The choice of $g$ can be
used to obtain information about the states $\mb{x} \in \Omega$ realized by the
stochastic trajectories. The choice of $f$ can be used to obtain information
about where the stochastic trajectories encounter the boundary $\partial
\Omega$.  The Dynkin formula and Backward-Kolmogorov equations provide
connections between these statistics and a collection of elliptic Partial
Differential Equations (PDEs) with boundary values problems of the form 
$\mathcal{L} u = -g, \,\, \mb{x} \in \Omega,\,\, u = f, \,\, \mb{x}
\in \partial \Omega$,~\cite{Oksendal2000}.}

\revI{We develop numerical methods for solving these PDEs on 
manifolds of general shape based on meshless methods.  
Meshless methods can be characterized broadly 
by their underlying discretizations.  This includes
Radial Basis Functions (RBF), Generalized Finite Differences (GFD), 
Moving Least Squares (MLS), and Reproducing Kernel
Particle Methods (RKPM),~\cite{gingold1977smoothed,buhmann2003radial,
lancaster1981surfaces,Liu1995,Kansa1990_1,Kansa1990_2}.  While
most meshfree approaches are for solutions of PDEs in euclidean
spaces $\mathbb{R}^d$ with $d=1,2,3$, recent work 
has focused on the manifold
setting, including~\cite{LowengrubLeung2011, liang2013solving,
suchde2019,Shankar_RBF_Lagrangian_2018, Mohammadi_Sphere_Kriging_LS_2019,
MacdonaldSurfPDEReactDiff2013,
ArroyoSurfMonge2019,PiretRBF_Orthog_Grad2012,FuselierWright2013}.
This includes methods for stablizing 
RBFs~\cite{fornberg2011stabilization,
flyer2012guide,fries2008convergence}, semi-lagrangian 
methods~\cite{Shankar_RBF_Lagrangian_2018}, approaches
avoiding use of surface 
coordinates~\cite{ShankarRBF_FD_RD_2015,
ShankarRBF_LOI_Aug_SurfPDE2018},
and methods that make use of the embedding 
space~\cite{PiretRBF_Orthog_Grad2012,PetrasClosestPointRBF2018,Ruuth2008}.  
This also includes methods using least-squares 
methods~\cite{Liang2013_SolversPointClouds,
GMLSStokes2020, TraskKuberry2020} and generalized finite difference
approaches~\cite{BertozziFourthOrderGeometries2006, Cheung2018_SurfPDE}.
}

\revI{For surfaces of general shape, 
we develop high-order numerical methods for computing statistics 
using Generalized Moving Least Squares (GMLS)
approximations~\cite{Mirzaei2012,wendland2004scattered}.  
Our approaches provide
meshless methods for solving PDEs on surfaces related
to~\cite{Liang2013_SolversPointClouds, GMLSStokes2020, TraskKuberry2020,
BertozziFourthOrderGeometries2006, Cheung2018_SurfPDE}. We focus here primarily 
on mean First Passage Times (FPTs), which are given by the special case $f = 0,\, g = 1$ with
$u(\mb{x}) = \mathbb{E}^{\mb{x}}\left[\tau\right]$.  Our methods also can be readily
extended for the cases with more general $f,g$.  We perform studies for a
variety of shapes to investigate how the methods converge.  We show our methods
can accurately capture both the surface geometry and the
action of the surface differential operators.  We then perform convergence
studies which show our methods have a high-order accuracy in approximating 
solutions of the surface PDEs.  
We show how our methods can be used to study 
how statistics are affected by the
surface geometry, drift dynamics, and spatially dependent diffusivities.  We also
discuss how the methods can be extended for computing more general 
path-dependent statistics for stochastic processes on surfaces.}

\section{First-Passage Times and Path-Dependent Statistics on Surfaces}

\begin{figure}[H]
\centering
\includegraphics[width=0.99\textwidth]{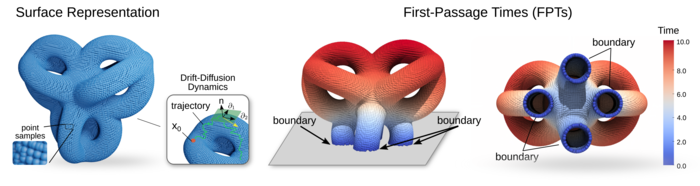} 
\hfill 
\caption{\revI{\textbf{Path-Dependent Statistics and Geometry.}
\textit{(left)} We develop meshless numerical methods for
solving surface PDEs associated with path-dependent statistics of
general drift-diffusion stochastic processes on surfaces of the form $d\mb{X}_t =
\mb{a}(\mb{X}_t)dt + \mb{b}(\mb{X}_t)d\mb{W}_t$.  
Statistics include the mean First Passage Times (FPTs) where the geometry and 
topology may play significant roles. \textit{(middle)} 
In FPTs, the process starts at $\mb{X}(0) = \mb{x}_0$ and the average 
time is computed for arriving at the boundary created by slicing 
the surface at $z = 0$. The surface geometry is represented by a collection
of point samples.  \textit{(right)} The surface colors indicate the value
of the FPT to reach the boundary when starting at a location on the surface.  
The results shown are for surface diffusion with 
$\mb{a}(\mb{x}) = 0$ and $\mb{b}(\mb{x}) = I$. }}
\label{fig:example1_ftp}
\end{figure}

The mean First-Passage Time (FPT) is defined as
\begin{eqnarray}
u(\mb{x}_0) = \mathbb{E}^{\small X(0)=\mb{x}_0}\left[\tau\right], \;\;
  \mbox{where } \tau = \inf\{t > 0 \, | \, \mb{X}_t \not \in \Omega\}.
\end{eqnarray}
This gives the time $\tau$ for a particle starting at location $X(0) = \mb{x}_0
\in \Omega$ to reach the boundary of the domain $\partial \Omega$.  We denote
by $u(\mb{x}_0)$ the average value of this passage time.  \revI{In
Figure~\ref{fig:example1_ftp}, we show the case of a curved surface with
boundary defined by a cut-plane at $z=0$.  The color at $\mb{x}$ indicates
$u(\mb{x})$ for the mean first passage time for a diffusion on the surface
starting at $\mb{x}$ to reach the boundary at $z = 0$, see
Figure~\ref{fig:example1_ftp}.}  

\revI{We consider the general drift-diffusion dynamics constrained to a surface
$\mathcal{M}$ governed by the Stochastic Differential Equations
(SDEs)
\begin{eqnarray}
\label{equ:sde}
d\mb{X}_t = \mb{a}(\mb{X}_t)dt + \mb{b}(\mb{X}_t) d\mb{W}_t, \hspace{0.5cm}
\mbox{subject to } \mb{X}_t \in \mathcal{M}. 
\end{eqnarray}
The $\mb{a}(\mb{x})$ term models the local drift in the dynamics and 
the $\mb{b}(\mb{x})$ the local diffusivity.  The $d\mb{W}_t$ denotes
increments of the Weiner process.
Throughout, we interpret our SDEs in terms of Ito stochastic
processes~\cite{Oksendal2000,Gardiner1985}. }
\lfoot{} 

\revI{While in principle stochastic simulations can be used to sample
trajectories for Monte-Carlo estimates of statistics, in practice
this can be expensive.  This arises from the need to 
reduce sufficiently the statistical sampling errors,
which can become particularly expensive when computing 
statistics for many initial starting locations.
Further challenges arise when
there are disparate time-scales in the dynamics resulting 
in numerical stiffness or in rare-events that result in slow decay of
the statistical errors~\cite{Newman1999,AtzbergerFPKMC2014}.  
We develop alternative approaches using Dynkin's 
formula and the Backward-Kolmogorov equations~\cite{Oksendal2000}.  This
relates the FPTs, and other path-dependent statistics, to elliptic PDE
boundary value problems.  We develop numerical methods for solving
these surface PDEs.}
  
For the stochastic process $\mb{X}_t$, we consider path-dependent statistics
$u$ of the form
\begin{eqnarray}
\label{equ:u_x_stat}
u(\mb{x}) = \mathbb{E}^{\mb{x}}\left[\int_0^\tau g(\mb{X}_t)dt \right]
+ \mathbb{E}^{\mb{x}}\left[f(\mb{X}_\tau)\right].
\end{eqnarray}
\revI{The $g(\mb{x})$ allows us to assign a value to individual trajectories based on the
locations they traverse inside the domain $\mb{x} \in \Omega$.
The $f(\mb{x})$ allows for assigning a
cost for the location where the process hits the boundary $\mb{x} \in \partial
\Omega$.}  The Dynkin formula for $q(\mb{x})$ is given by~\cite{Oksendal2000} 
\begin{eqnarray}
\label{subsec_DynkinFormula}
  \mathbb{E}^{\mb{x}}\left[q(\mb{X}_\tau)\right] = q(\mb{x}) 
  + \mathbb{E}^{\mb{x}}\left[\int_0^\tau \mathcal{L}q(\mb{X}_s)ds \right].
\end{eqnarray}
\revI{The $\mathcal{L}$ is the 
Infinitesimal Generator associated with the SDE in
equation~\ref{equ:sde}.  This can be expressed as the differential operator}
\begin{eqnarray}
\label{equ:inf_gen}
\mathcal{L} =  \mb{a}\cdot \frac{\partial}{\partial \mb{x}}
+ \frac{1}{2} \mb{b}\mb{b}^T :\frac{\partial^2}{\partial \mb{x}^2}.
\end{eqnarray}
\revI{The $\cdot$ denotes the usual dot product with a gradient, and the $:$
denotes the Hadamard element-wise dot product (tensor contraction over
indices).}  We take $q(\mb{x}) = f(\mb{x}),\; 
\mb{x}\in \partial \Omega$ and $\mathcal{L}q(\mb{x}) =
-g(\mb{x}), \; \mb{x} \in \Omega$. After rearranging terms, we have formally $u(\mb{x}) =
q(\mb{x}), \; \mb{x} \in \bar{\Omega}$.  \revI{Using this and equation~\ref{equ:inf_gen}, we
can express the statistic $u$ as the solution of a surface elliptic PDE.  This gives
the PDE boundary-value problem of the form}
\begin{eqnarray}
\label{equ_surf_PDE_bk}
\mathcal{L} u = -g, \,\, \mb{x} \in \Omega\,\,\;\; u = f, \,\, \mb{x} \in \partial \Omega,
\end{eqnarray}
where $\mathcal{L}$ is given in equation~\ref{equ:inf_gen}.  We require throughout $f \in C_0^2(\partial \Omega)$ and $g \in C_0^2(\Omega)$.

\revI{In the case with $f = 0$ and $g = -1$, we obtain the mean First-Passage Times (FPTs) 
\begin{eqnarray}
u(\mb{x})  = \mathbb{E}^{\mb{x}}\left[\tau\right],\;\;
\tau = \inf \{t > 0 | \mb{X}_t \neq \Omega\}.
\end{eqnarray}
This is the amount of time on average it takes for the drift-diffusion process starting at location $\mb{x}$ to hit the boundary of the domain $\partial \Omega$.}

\revI{To compute these and other path-dependent statistics for stochastic processes
constrained to surfaces, we develop solvers for the surface PDEs of
equation~\ref{equ_surf_PDE_bk}.  We show how these methods can be used to
investigate the roles played by geometry and the drift-diffusion dynamics
in the statistics of equation~\ref{equ:u_x_stat}.}

\subsection{Stochastic Processes and Drift-Diffusion Dynamics on Surfaces}
\label{sec:dynamics_ab}
\revI{We formulate descriptions for Brownian motion and other 
more general stochastic processes on surfaces. 
We express the surface drift-diffusion dynamics using
the notation
$d\mb{X}_t = \mb{a}dt + \mb{b}d\mb{W}_t.$
The $\mb{a} = \mb{a}(\mb{x}) \in \mathbb{R}^n$ and 
$\mb{b} = \mb{b}(\mb{x}) \in
\mathbb{R}^{n \times m}$.  The $\mb{X}_t \in \mathbb{R}^n$ defines 
the location in terms of the ambient embedding
space. }  In the embedding space we assume the drift $\mb{a}$
is in the tangent plane at each point $\mb{a}(\mb{x}) \in
T\mathcal{M}_{\mb{x}}$ and that the range of the diffusion tensor $\mb{b}$ is
in the tangent space $\mbox{range}(\mb{b}) \subset T\mathcal{M}_{\mb{x}}$. \revI{As
a result, the null-space of $\mb{b}^T$ includes the space of vectors orthogonal to
the tangent space $T\mathcal{M}_{\mb{x}}$.  This has the consequence that
$\mb{b}\mb{b}^T$ has a null-space of vectors orthogonal to the tangent space.
This also has the range of vectors that are at most spanned by $\mb{b}$ in the 
tangent space, see Figure~\ref{fig:example1_ftp}.}

In modeling systems, it is often convenient to specify the drift and diffusion
using the embedding space representations of $\mb{a}$ and $\mb{b}$.  \revI{In
numerical methods, and when performing other practical calculations, it is often
convenient to express the dynamics using local coordinate charts $\mb{q} =
(q_1,q_2)$ with embedding map $\mb{x} = \bsy{\sigma}(q_1,q_2)$.}  It will be
convenient to have ways to convert between these types of descriptions.  
\revI{In local coordinate charts, the drift-diffusion dynamics 
can be expressed as}
\begin{eqnarray}
d\mb{q}_t = \bsy{\alpha}dt + \bsy{\beta}d\tilde{\mb{W}}_t,
\end{eqnarray}
where $\bsy{\alpha}(\mb{q}) \in \mathbb{R}^2$ and $\bsy{\beta}(\mb{q}) \in
\mathbb{R}^{2 \times \tilde{m}}$ with $\tilde{m}$ the number of noise sources.

These descriptions are connected through the embedding map $\bsy{\sigma}$ given
by $\mb{X} = \bsy{\sigma}(\mb{q})$.  \revI{We use the Ito Lemma~\cite{Oksendal2000}
to connect the dynamics in the embedding space with those in the local 
coordinate chart to obtain}
\begin{eqnarray}
d\mb{X}_t & = & \left[\alpha^a + \frac{1}{2}\Gamma_{ij}^a \bsy{\beta}^i  \bsy{\beta}^{j,T} 
\right] \bsy{\sigma}_{q^a}  dt 
+ 
\bsy{\sigma}_{q^a} \bsy{\beta}^a d\tilde{\mb{W}}_t .
\end{eqnarray}
\revI{The $\Gamma_{ij}^a$ are the Christoffel symbols of the 
surface~\cite{Pressley2001,Abraham1988}.  We take
$\bsy{\beta}^a \in \mathbb{R}^2$ to be a row vector and 
$d\tilde{\mb{W}}_t \in \mathbb{R}^2$ to be a column vector.  
When these combine they give a scalar for each of the 
coordinate components.  We also use the notational conventions 
that $\partial_{q^a} = \partial
\bsy{\sigma} / \partial q^a = \bsy{\sigma}_{q^a}$ and
that repeated indices sum.  We can use $\{\partial_{q^a}\}_a$ 
to represent the dynamics of $\mb{X}_t$ in the 
surrounding embedding space.  We discuss further details on
how to compute these terms from the surface geometry in
Appendix~\ref{appendix:monge_gauge_param}.}

We obtain relations between descriptions using the expressions for the dynamics
in the embedding space $\mathbb{R}^3$ and the dynamics expressed in the
coordinates $\mb{q}$. \revI{These can be summarized as 
\begin{eqnarray}
\mb{b}\mb{b}^T = \bsy{\beta}^a\bsy{\beta}^{b,T}
\partial_{q^a} \partial_{q^{b}}^T = \bsy{\beta}^a\bsy{\beta}^{b,T}
\bsy{\sigma}_{q^a} \bsy{\sigma}_{q^{b}}^T, 
\hspace{0.6cm}
a^c  = g^{cb}\langle \mb{a}, \bsy{\sigma}_{q^{b}} \rangle_g, 
\hspace{0.6cm}
\alpha^c =  a^c - \frac{1}{2} \Gamma_{ab}^c
\bsy{\beta}^a\bsy{\beta}^{b,T}. \label{eq_alpha_c_coeficient}
\end{eqnarray}
The $\bsy{\beta}^a\bsy{\beta}^{b,T}$ is a scalar, since $\bsy{\beta}^c$ is a
row vector.}  The $\langle \cdot, \cdot \rangle_g$ denotes the inner-product in
the embedding space.  The $a^c$ is useful so that we can represent the vector
as $\mb{a} = a^c \bsy{\sigma}_{q^c}$.  We still subscript with $g$ even in the
embedding space to highlight that this corresponds to the metric inner-product
when using a coordinate chart.  In practice, this is computed readily if we
represent the vectors in the ambient space as the usual Euclidean inner-product
between vectors.  

\revI{While there are many possible choices of the noise terms
$\bsy{\beta}$ and $\mb{b}$ consistent with these equations, each of them
are weakly equivalent since they have the same marginal probability densities.}  For
convenience, we make the specific choice 
\begin{eqnarray}
\mb{b} = \bsy{\sigma}_{q^a} \bsy{\beta}^a.
\end{eqnarray}
This can be verified to satisfy the relations above.  

The $\sigma_{q^a}$ is a column vector and $\bsy{\beta}^a$ is a row vector so
that $\mb{b}$ is a $2\times 2$ matrix.  We choose also for our driving Brownian
motion here $W_t = \tilde{W}_t$ with $m = \tilde{m} = 2$.  With this
correspondence, we also can express $\bsy{\beta}^a$ in terms of $\mb{b}$ as
\begin{eqnarray}
\bsy{\sigma}_{q^b}^T \mb{b} = \bsy{\sigma}_{q^b}^T \bsy{\sigma}_{q^a} \bsy{\beta}^a = g_{ab} \bsy{\beta}^a.
\end{eqnarray}
This can be solved using the inverse metric tensor as
\begin{eqnarray}
\bsy{\beta}^c = g^{cb}\bsy{\sigma}_{q^b}^T \mb{b}.
\end{eqnarray}
The $\bsy{\sigma}_{q^b}^T$ is a row vector, $\mb{b}$ is a matrix, and this
yields the row vector $\bsy{\beta}^c$.  

We can substitute this in equation~\ref{eq_alpha_c_coeficient} to relate the
coordinate chart drift and diffusion terms to the embedding space quantities as 
\begin{eqnarray}
\label{equ:ambient_to_chart}
\bsy{\beta}^c  =  g^{cb}\bsy{\sigma}_{q^b}^T \mb{b}, \;\;
a^c =  g^{cb}\langle \mb{a}, \partial_{q^b} \rangle_g, \;\;
\alpha^c  =  a^c - \frac{1}{2} \Gamma_{ab}^c \bsy{\beta}^a\bsy{\beta}^{b,T}.
\end{eqnarray}
\revI{The 
Infinitesimal
Generator $\mathcal{L}$ for the surface drift-diffusion
process in equation~\ref{equ:inf_gen} can be expressed in the local coordinate
chart as}
\begin{eqnarray}
\label{equ_dynkin_expanded}
\mathcal{L} u = \alpha^c\frac{\partial u}{\partial q^c}
+ \frac{1}{2} \bsy{\beta}^{c}\bsy{\beta}^{d,T} \frac{\partial^2 u}{\partial q^c\partial q^d}.
\end{eqnarray}
In practice, this is used to compute the action of $\mathcal{L}$ on $u$ when
evaluating terms in the surface PDEs in equation~\ref{equ_surf_PDE_bk}.  

We remark that a convenient feature of these calculations is that in the final
expressions no derivatives of the drift and diffusion coefficients were needed.
This is particularly helpful on complicated surfaces for practical
calculations.  For modeling and simulation in practice, we use the convention
for the data structures that $\mb{a}$ is given as a vector field on the point
cloud.  We also specify that $\mb{b}$ is given in terms of components
$\mb{b}[1]$,$\mb{b}[2]$,$\mb{b}[3]$, where $\mb{b}[i]$ is the $i^{th}$ column
of the matrix.  In this way, we can represent readily the tensors
$\mb{a},\mb{b}$ and terms needed in the calculations of $\mathcal{L}$ using
equations~\ref{equ:ambient_to_chart} and~\ref{equ_dynkin_expanded} when solving
equation~\ref{equ_surf_PDE_bk}.

\section{GMLS Solvers for Surface Partial Differential Equations (PDEs)}
\label{sec:gmls_solvers}
To obtain first passage times and other path-dependent statistics, we need
numerical methods for solving the elliptic surface PDEs in
equation~\ref{equ_surf_PDE_bk}.  This requires approximation of the surface
geometry and associated differential operators arising in
equations~\ref{equ_surf_PDE_bk},~\ref{equ:ambient_to_chart},
and~\ref{equ_dynkin_expanded}.  This poses challenges given the high-order
derivatives required of geometric quantities such as the local surface
curvature and metric tensor.  We address this by developing meshless approaches
based on Generalized Moving Least Squares
(GMLS)~\cite{wendland2004scattered,Mirzaei2012}.  \revI{We use this to obtain
high-order approximations for geometric quantities associated with the local
surface geometry.}  
\revI{In
meshless methods, a collection of point samples is used to represent the
manifold geometry and at these locations the surface scalar and vector 
fields, see
Figures~\ref{fig:example1_ftp} and~\ref{fig:mls}.}  
\revI{We develop GMLS approaches for estimating both the local
geometry and surface differential operators which are used to 
build collocation methods for the surface PDEs.}

\subsection{Generalized Moving Least Squares (GMLS) Approximation}
\label{sec_gmls}

The manifold is represented as a point cloud $\{x_i\}_{i=1}^N$ from which we
need to approximate associated differential operators and geometric quantities.
Generalized Moving Least Squares (GMLS) approximates operators of the
underlying surface fields by solving a collection of local least-squares
problems~\cite{wendland2004scattered,Mirzaei2012}.  \revI{Given at
each $\mb{x}_i$ a finite dimensional Banach space $\mathbb{V}$ and dual space 
$\mathbb{V}^*$, an approximation in $\mathbb{V}^*$ is sought for a 
target operator $\tau_{\mb{x}_i}[u]$.
The $\tau_{\mb{x}_i} \in \mathbb{V}^*$ and $u(\mb{x}) \in \mathbb{V}$ for $\mb{x} \in \Omega
\subset \mathbb{R}^d$ where $\Omega$ is a compact domain.  In practice, we take
$\mathbb{V} = \mathbb{V}_n$ to be the collection of multinomials up 
to total degree $n$.
To relate $u \in \mathbb{V}$ to a representative $p^* \in \mathbb{V}_n$, we
consider a collection of probing linear functionals $\Lambda =
\{\lambda_j\}_{j=1}^N$ that serve to characterize $u$ by $\Lambda[u] =
(\lambda_1[u],\lambda_2[u],\ldots,\lambda_N[u])$.}  We construct the
approximation $p^*$ by solving the $\ell^2$-optimization problem
\begin{eqnarray} p^* = \arg\min_{p \in \mathbb{V}_n} \sum_{j = 1}^N
\left(\lambda_j[u] - \lambda_j[p]\right)^2 \omega(\lambda_j,\tau_{x_i}).
\label{eqn_p_star} \end{eqnarray}
The $\omega(\lambda_j,\tau_{x_i})$ provides
weights characterizing the importance of the particular sampling function
$\lambda_j$ in estimating $\tau_{x_j}$.  For example, we can take $\lambda_j =
\delta(x - x_j)$ with $\lambda_j[u] = u(x_j)$ and we can take
$\omega(\lambda_j,\tau_{x_i}) = \omega(x_j - x_i)$ to be a function that decays
as the distance increases between $x_j$ and $x_i$.  We can further take $\omega
= \omega(\|x_j - x_i\|)$ to be a radial function that decays to zero when $r >
\epsilon$.  In practice, we use $\omega(\|x_j - x_i\|) = \left(1 - \|x_j -
x_i\|/\epsilon\right)_+^p$ where $\left(\cdot\right)_+ = \max(\cdot,0)$.

\revI{Consider the basis functions $\Phi = \{\phi_i\}_{i = 1}^{d_n}$ for
$\mathbb{V}_n$ so that $\mathbb{V}_n =
\mbox{span}\{\phi_1,\phi_2,\ldots,\phi_{d_n}\}$ with $d_n = \dim \mathbb{V}_n$.}
Any function $p \in \mathbb{V}_n$ now can be expressed as
\begin{eqnarray}
p(\mb{x}) = \Phi^T \mb{a} = \sum_{i=1}^{d_n} a_i \phi_i(\mb{x}).
\end{eqnarray}
Let $\tau[\Phi]$ denote a vector where the components are the target operator
acting on each of the basis elements.  For the the target operator $\tau$, we
obtain the GMLS approximation $\tilde{\tau}$ by considering how $\tau$ acts on
the optimal reconstruction $p^* = \Phi^T \mb{a}^*$ of $u$,
\begin{eqnarray}
\label{equ:gmls_tau}
\tilde{\tau}[u] := \tau[p^*] = \tau[\Phi]^T \mb{a}^*.
\end{eqnarray}
Conditions ensuring the existence of solutions to equation \ref{eqn_p_star}
depend primarily on the unisolvency of $\Lambda$ over $\mathbb{V}$ and
distribution of the point samples $\{\mb{x}_i\}$.  For theoretical results related
to GMLS see~\cite{Mirzaei2012,wendland2004scattered}.  GMLS has been primarily
used to obtain approximations of constant coefficient linear differential
operators~\cite{wendland2004scattered}.

\begin{figure}[H]
\centering
\includegraphics[width=0.99\textwidth]{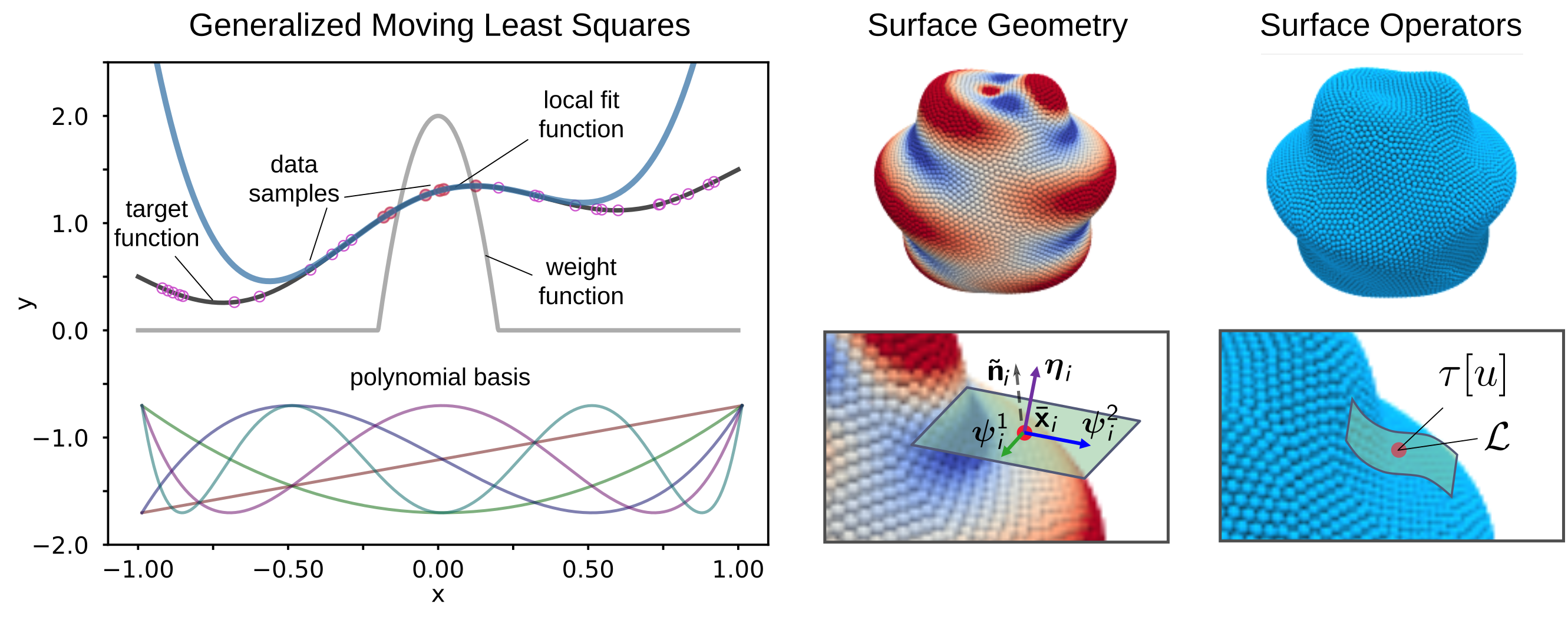} 
\hfill 
\caption{\revI{\textbf{Generalized Moving Least Squares (GMLS) on Surfaces.}
\textit{(left)} We approximate the geometry and the
action of operators $\tau[u]$ on general surfaces using 
Generalized Moving Least Squares (GMLS).
  To approximate a target function $u(\mb{x})$ in the neighborhood of a given point
  $\mb{x}^*$, GMLS uses the data samples
$\{(\mb{x}_j,u(\mb{x}_j))\}_{j=1}^m$, a collection of probing 
functionals $\{\lambda_k[u]\}$, and weights to localize the fit
$w(\mb{x}^*,\mb{x}_j)$.  These are used to fit locally a function
$p^* = p^*(\mb{x}) = p^*(\mb{x};\mb{x}^*) \in \mathbb{V}$.  
Multivariate polynomials of finite degree are used for 
the function space $\mathbb{V}$.   
\textit{(middle)}
Operators are approximated by locally solving related
least-squares problems using equations~\ref{eqn_p_star} 
and~\ref{equ:gmls_tau}.  To handle the geometry at $\mb{x}^*$,
Principle Component Analysis (PCA) is performed using local points
to obtain an approximate tangent plane.  Geometric quantities
are estimated by using GMLS to fit locally the surface
in the Monge-Gauge and approximating the related differential
operators.  Shown is the GMLS estimate of the Gaussian curvature
on the surface.
\textit{(right)}  The geometric results are used further
in conjunction with the GMLS approach to estimate 
the action of the surface 
differential operators $\tau[u] = \mathcal{L}u$.  
Expressions for obtaining geometric quantities in terms of the 
local surface coordinates are given in 
Appendix~\ref{appendix:monge_gauge_param}.}}
\label{fig:mls}
\end{figure}
For our first-passage time problems, the target operators are technically
non-linear given their dependence on the geometry of the underlying manifold
which also must be estimated from the point cloud.  We handle this using a two
stage approximation approach.  In this first stage, we use the point cloud to
estimate two basis vectors $\psi_1$ and $\psi_2$ for the local tangent plane
using Principle Component Analysis
(PCA)~\cite{Pearson_PCA_1901,HastieElementsStatLearning2001}.  We use these
basis vectors to construct a local coordinate chart $(q_1,q_2)$. In this chart,
we fit a function $p^* \in \mathbb{V}$ to the local point cloud to obtain a
Monge-Gauge~\cite{Pressley2001} representation of the surface
$(q_1,q_2,p^*(q_1,q_2))$.  We use GMLS to estimate target geometric quantities
of interest, such as the local Gaussian curvature or high-order derivatives.
In the second stage, we use the estimated geometric quantities to specify the
target operators $\tau[u]$ for the fields $u$.  We then use again GMLS to
obtain an approximation of $\tilde{\tau}$ and to compute numerically the target
operator values at $\mb{x}_i$.  We have used related procedures for solving
hydrodynamic equations on manifolds in~\cite{GMLSStokes2020}.  \revI{We illustrate
this approximation approach in Figure~\ref{fig:mls}.}


\subsection{Numerical Methods for Solving the Surface PDEs}

To compute numerically the first passage times, we develop solvers for the elliptic PDE boundary value problems given in equation~\ref{equ_surf_PDE_bk}.  This is organized by representing the $f$ and $g$ function inputs and the $\mb{a}$ and $\mb{b}$ fields specifying the drift-diffusion dynamics.  To avoid complications with local coordinate charts, we numerically represent all input data globally using the ambient embedding space coordinates $\mb{x} \in \mathbb{R}^3$.  The $\mb{a} \in \mathbb{R}^3$ with only the tangential components playing a role in practice.  Similarly, the $\mb{b}$ tensor is represented by three vector fields $\mb{b_1}$, $\mb{b_1}$, and $\mb{b_3}$.  We use labels on the point cloud to determine which regions are to be considered interior to $\Omega$ and which are part of the boundary of $\partial \Omega$.  We only require $f$ to be evaluated on $\partial \Omega$, while $g$ must return reliable values for all $\mb{x} \in \Omega$.

We use our GMLS methods in section~\ref{sec_gmls} to estimate the surface geometric quantities and the action of the operator in equation~\ref{equ_surf_PDE_bk}.  This allows us to construct at each $\mb{x}_i$ an equation for relating $\mathcal{L}u(\mb{x}_i) = -g(\mb{x}_i)$.   This provides our collocation method for determining $u(\mb{x}_i)$.  
Let $[\tilde{u}]_i = u(\mb{x}_i)$ and $[\tilde{g}]_i = g(\mb{x}_i)$.
Collecting these equations together gives a sparse linear system 
$
A \tilde{\mb{u}} = \tilde{\mb{g}}.
$
We solve this large sparse linear system using GMRES with algebraic multigrid (AMG) preconditioning using Trilinos~\cite{trilinos2005}.  Our solvers have been implemented within a framework for GMLS problems using the Compadre library and PyCompadre~\cite{Compadre2019}.  The toolkit provides domain decomposed distributed vector representation of fields as well as global matrix assembly.
The capability of the library were also extended for our surface geometry calculations by implementing symbolically generated target operators.  Our methods also made use of the iterative block solvers of (Belos~\cite{Belos}), block preconditioners of (Teko) and the AMG preconditioners of (MueLu~\cite{MueLu2014a,MueLu2014b}) within the Trilinos software framework~\cite{trilinos2005}.  The framework facilitates developing a scalable implementation of our methods providing ways to use sparse data structures, parallelization, and hardware accelerations.


\section{Convergence Results}

We investigate the convergence of our GMLS solvers developed in section~\ref{sec:gmls_solvers}.  The target operators that arise in the surface PDE boundary-value problems involve a non-linear approximation.  This arises from the coupling between the contributions to the error from the differential terms of the surface operators and the GMLS estimations used for the surface geometry.

\subsection{Surface Geometries for Validation Studies}
\begin{figure}[H]
\centerline{\includegraphics[width=0.99\columnwidth]{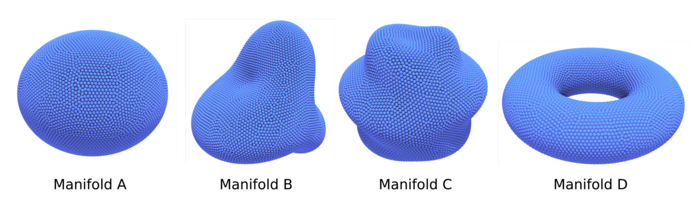}}
  \caption{\revI{\textbf{Point Sample Representations of the Surfaces.}} }
\label{fig:manifolds}	
\end{figure}

We perform studies using four different surface geometries: (i) ellipsoid, (ii)
radial manifold I, (iii) radial manifold II, and (iv) torus.  We label these as
Manifold A--D, see Figure~\ref{fig:manifolds}.    We study convergence as the
surface point sampling is refined.  We take each refinement to have
approximately four times the number of points as the previous level.  This aims
to have the fill distance $h$ halve under each refinement.

The manifolds can be described by the following implicit equations.  Manifold A
is an ellipsoid defined by the equation ${x^2}/{a^2} + {y^2}/{b^2} + z^2 =
s_0^2$ with $a = 1.2, b = 1.2, s_0^2 = 1$.  Manifold B is a radial manifold
defined in spherical coordinates by $(\theta,\phi,r(\theta,\phi)$ where
$r(\theta,\phi) = 1 + r_0\sin(3\phi)\cos(\theta)$ with $r_0 = 0.1$.  Manifold C
is a radial manifold defined in spherical coordinates by
$(\theta,\phi,r(\theta,\phi)$ where $r(\theta,\phi) = 1 +
r_0\sin(7\phi)\cos(\theta)$ with $r_0 = 0.1$.  Manifold D is a torus defined by
the equation $(s_1^2 - \sqrt{x^2 + y^2} )^2 + z^2 = s_2^2$ with $s_1^2 = 0.7,
s_2^2 = 0.3$.  Each of the manifolds shown are represented by quasi-uniform
point sets with approximately $n = 10^4$ samples.  For quasi-uniform sampling
we expect the fill-distance $h$ to scale as $h \sim 1/\sqrt{n}$.  We report our
results throughout using the notation $\bar{h}^{-1} = \sqrt{n}$. Additional
information on the number $n$ used in samplings can be found in
Appendix~\ref{appendix:sampling_res_manifolds}.

In our first passage time calculations, unless indicated otherwise, we generate
the boundary $\partial \Omega$ as the points at $z = 0$.  In practice, we treat
any points with $|z| < 4 \times 10^{-2}$ as boundary points. This serves to
thicken the boundary region and at points $\mb{x}_i$ near the boundary
$\partial \Omega$ helps achieve unisolvency in the local least-squares
problems.

\subsection{Test Fields and Manufactured Solutions}
\label{sec:test_u_ab}
We study the accuracy of the our GMLS approximations of the surface operators
by investigating their action on test fields $\hat{\mb{u}}(\mb{x})$.  We
generate $u$ by using smooth functions $\hat{\mb{u}}(x,y,z)$ parameterized in
the embedding space $\mathbb{R}^3$.  By the smoothness of the surface manifolds
$\mathcal{M}$, evaluation at $\mb{x} \in \mathcal{M}$ gives a smooth surface
field $u$.  More formally, this corresponds to using the inclusion map $\iota:
\mathbb{R}^3 \xhookrightarrow{} \mathcal{M}$ to obtain $\mb{u}(\mb{x}) =
\iota_{\mb{x}} \hat{\mb{u}}(\cdot)$.   We also use this inclusion map approach
to generate test surface drift tensors $\mb{a}$ and diffusion tensors $\mb{b}$
for the stochastic process $\mb{X}_t$ arising in equation~\ref{equ:sde}
and~\ref{equ:inf_gen}.

For our validation studies, we generate test fields $u$ for the Manifolds A--D
using $\mb{u}(\mb{x}) = z(x^4 + y^4 - 6x^2y^2)$, which is an extension of a
degree $5$ Spherical Harmonic to $\mathbb{R}^3$.  For the drift-diffusion
tensors for Manifold A, we use \begin{eqnarray} \mb{a}_A = y
  \tilde{\bs{\sigma}}_1 + xz \tilde{\bs{\sigma}}_1, \;\;\; \mb{b}_A =
  \begin{bmatrix} xz \tilde{\bs{\sigma}}_1 + y^2 \tilde{\bs{\sigma}}_2, xyz
  \tilde{\bs{\sigma}}_1 + \cos(y^2) \tilde{\bs{\sigma}}_2, z
  \tilde{\bs{\sigma}}_1 + (y^3 + x) \tilde{\bs{\sigma}}_2 \end{bmatrix} .
\end{eqnarray} The $\tilde{\bs{\sigma}}_1, \tilde{\bs{\sigma}}_2$ provide a
local orthogonal tangent basis at every point on the Manifold A.  In practice
in our numerical calculations, this does not need to be smooth in $\mb{x}$ and
we construct these as convenient using our local tangent plane approximation
and Gram-Schmidt orthogonalization~\cite{Burden2010}.  

For the Manifolds B and C, we use for the drift-diffusion tensors tangential
projections of the vector fields in $\mathbb{R}^3$ given by
\begin{eqnarray}
\mb{a}_{B,C} = \begin{bmatrix} y,\; xz,\; x^2 yz \end{bmatrix}^T,\,\,\,\,\,\,
\mb{b}_{B,C}  =  \begin{bmatrix}(xz, xyz, z), \;
                        (y^2, \cos(y^2), y^3+x), \;
                        (x^2+ y, e^{-z}, e^{y}) 
                        \end{bmatrix}^T.
\end{eqnarray}
For Manifold D, we use drift-diffusion tensors
\begin{eqnarray}
\mb{a}_{D} = y \bs{\sigma}_u + xz \bs{\sigma}_v,\;\;\; \mb{b}_{D} =
  \begin{bmatrix} \sin(xy)\bs{\sigma}_u + y^2  \bs{\sigma}_v , xyz\bs{\sigma}_u
  + \cos(y^2) \bs{\sigma}_v , e^{z}\bs{\sigma}_u + (y^3 +x)\bs{\sigma}_v
  \end{bmatrix}.
\end{eqnarray}
The $\bs{\sigma}_u, \bs{\sigma}_v$ provide a basis for the tangent space smooth
in $\mb{x}$ given by the directional derivatives of the global parameterization
$\bs{\sigma}(u,v)$ of the torus~\cite{Pressley2001}.  We characterize the
accuracy of the GMLS approximation $\tilde{L}$ of the surface operators
$\mathcal{L}$ using the $\ell^2$-error \begin{eqnarray}
  \epsilon_{op}^2\left(\tilde{L}\right) &=& \left\|\tilde{L} u - \mathcal{L} u
  \right\|_{n,2}^2 = \frac{1}{n} \sum_{i=1}^n \left(\tilde{L} u(\mb{x}_i) -
  \mathcal{L} u(\mb{x}_i)\right)^2.  \end{eqnarray} In these studies, we
  evaluate to high precision the action of the operators $\mathcal{L}$ by
  symbolic calculations using SymPy~\cite{Sympy2017}.  In general, we emphasize
  that such calculations of expressions symbolically can be prohibitive. \revI{Using
  this approach, we investigate the accuracy of the GMLS approximation of the
  operator for each of the manifolds.  We use the approximation spaces with
  multivariate polynomials of degree $m \in \{2, 4, 6\}$.}  These results are
  reported in Tables~\ref{table:DynkinPequals2}--~\ref{table:DynkinPequals6}.

\subsection{Results of Convergence Studies}
We performed convergence studies for the Manifolds A--D shapes in
Figure~\ref{fig:manifolds} with the test fields $u$ and drift-diffusion tensors
$\mb{a}, \mb{b}$ discussed in Section~\ref{sec:test_u_ab}.  \revI{As the manifold
resolution is refined, we study the convergence of the GMLS approximations of
the surface operators.  For the approximation spaces $\mathbb{V}$,
we use multivariate polynomials of
degree $m \in \{2,4,6\}$.}  We estimate the convergence rates using the log-log
slope of the error $\epsilon_{op}$ as the fill distance parameter $h$ is varied
between levels of refinement.  We report these results in
Tables~\ref{table:DynkinPequals2}--~\ref{table:DynkinPequals6}.

\revI{We find in our empirical results that when using GMLS with multinomial spaces of
degree $m$ to evaluate the elliptic PDEs operator $\mathcal{L}$ of order $k$,
we obtain convergence results of order $m$.}  \revI{We remark that the
standard theory does not apply for our surface operators, 
since there are non-linearities from the coupling
between the surface geometry and the differential terms in the operator.}  
Our results
are suggestive that our surface operator approximations achieve results
exhibiting a trend consistent with simpler differential operators approximated
by GMLS~\cite{Mirzaei2012}.  This suggests our methods are resolving the
surface geometric contributions with sufficient precision to minimize their
contributions to the error.  Since the operator $\mathcal{L}$ is of
differential order $k=2$ we use throughout approximating polynomial spaces of
at least degree $m \geq k \geq 2$.

\begin{table}[H]
\label{DynkinPequals2}
\begin{center}
\begin{tabular}{ccccccccccc}
\textbf{}  & \multicolumn{2}{c}{\textbf{Manifold A}}      & \multicolumn{2}{c}{\textbf{Manifold B}} & \multicolumn{2}{c}{\textbf{Manifold C}} & & & \multicolumn{2}{c}{\textbf{Manifold D}}      \\\cline{2-7}\cline{10-11}
\textbf{h} & \textbf{$\ell_2$-error} & \textbf{Rate} & \textbf{$\ell_2$-error} & \textbf{Rate} & \textbf{$\ell_2$-error} & \textbf{Rate} & & \textbf{h} & \textbf{$\ell_2$-error} & \textbf{Rate} \\
0.1 & 2.6503e-02 & - & 3.6072e-02 & -  & 1.1513e-02  & -&  & .08 &  8.1920e-03 & -  \\
0.05 & 5.8961e-03 & 2.14 & 8.3322e-03 & 2.11 & 2.3778e-03 & 2.27 & & .04 & 1.8071e-03 & 2.21 \\
0.025 & 1.4951e-03 & 1.97 & 2.1451e-03 & 1.95 & 6.3865e-04 & 1.90 & & .02 &  4.3928e-04 & 2.06 \\
0.0125 & 3.5399e-04 & 2.07 & 5.1159e-04 & 2.06 & 1.6225e-04 &  1.98 & & .01 & 1.0461e-04 & 2.05 \\
\end{tabular}
\end{center}
\caption{GMLS Convergence Rates for $\mathcal{L}$.  For our GMLS methods using approximation spaces based on multivariate polynomials of degree $m = 2$, we find convergence rate of $\sim 2^{nd} \mbox{order}$.}
\label{table:DynkinPequals2}
\end{table}

\begin{table}[H]
\label{DynkinPequals4}
\begin{center}
\begin{tabular}{ccccccccccc}
\textbf{}  & \multicolumn{2}{c}{\textbf{Manifold A}}      & \multicolumn{2}{c}{\textbf{Manifold B}} & \multicolumn{2}{c}{\textbf{Manifold C}} & & & \multicolumn{2}{c}{\textbf{Manifold D}}      \\\cline{2-7}\cline{10-11}
\textbf{h} & \textbf{$\ell_2$-error} & \textbf{Rate} & \textbf{$\ell_2$-error} & \textbf{Rate} & \textbf{$\ell_2$-error} & \textbf{Rate} & & \textbf{h} & \textbf{$\ell_2$-error} & \textbf{Rate} \\
0.1 & 5.8394e-03 & - & 2.5843e-02 & -  & 5.5971e-03  & -&  & .08 &  2.5936e-03 & -  \\
0.05 & 2.0355e-04 & 4.78 & 1.1594e-03 & 4.48 & 3.7890e-04 & 3.88 & & .04 & 1.0359e-04 & 4.72 \\
0.025 & 1.1382e-05 & 4.14 & 8.2945e-05 & 3.80 & 1.9628e-05 & 4.27 & & .02 &  5.4028e-06 & 4.30 \\
0.0125 & 6.6746e-07 & 4.08 & 5.1207e-06 & 4.01 & 1.4226e-06 &  3.79 & & .01 & 3.0733e-07 & 4.11 \\
\end{tabular}
\end{center}
\caption{GMLS Convergence Rates for $\mathcal{L}$.  For our GMLS methods using approximation spaces based on multivariate polynomials of degree $m = 4$, we find convergence rate of $\sim 4^{th} \mbox{order}$.}
\label{table:DynkinPequals4}
\end{table}

\begin{table}[H]
\label{DynkinPequals6}
\begin{center}
\begin{tabular}{ccccccccccc}
\textbf{}  & \multicolumn{2}{c}{\textbf{Manifold A}}      & \multicolumn{2}{c}{\textbf{Manifold B}} & \multicolumn{2}{c}{\textbf{Manifold C}} & & & \multicolumn{2}{c}{\textbf{Manifold D}}      \\\cline{2-7}\cline{10-11}
\textbf{h} & \textbf{$\ell_2$-error} & \textbf{Rate} & \textbf{$\ell_2$-error} & \textbf{Rate} & \textbf{$\ell_2$-error} & \textbf{Rate} & & \textbf{h} & \textbf{$\ell_2$-error} & \textbf{Rate} \\
0.1 & 1.8032e-05 & - & 2.0087e-03 & -  & 2.1854e-03  & -&  & .08 &  2.7306e-04 & -  \\
0.05 & 2.3775e-07 & 6.16 & 5.4584e-05 & 5.21 & 1.1625e-04 & 4.24 & & .04 & 2.6141e-06 & 6.81 \\
0.025 & 4.3492e-09 & 5.75 & 1.4770e-06 & 5.20 & 1.5943e-06 & 6.20 & & .02 &  7.4551e-08 & 5.18 \\
0.0125 & 6.4443e-11 & 6.06 & 1.9655e-08 & 6.22 & 3.1347e-08 &  5.68 & & .01 & 1.6423e-09 & 5.47 \\
\end{tabular}
\end{center}
\caption{GMLS Convergence Rates for $\mathcal{L}$.  For our GMLS methods using approximation spaces based on multivariate polynomials of degree $m = 6$, we find convergence rate of $\sim 6^{th} \mbox{order}$.}
\label{table:DynkinPequals6}
\end{table}


\section{Results for First Passage Times on Surfaces}
For first passage times on surfaces, we investigate the roles played by
the geometry, drift dynamics, and spatial-dependence of diffusivity.  For our studies, we consider the general surface Langevin dynamics 
\begin{eqnarray}
\label{equ:langevin_surf}
d\mb{X}_t =  -\frac{1}{\gamma} \gradX U(\mb{X})dt + \sqrt{2D} \mb{I} d\mb{W}_t \\
\nonumber
\mbox{subject to } \mb{X}_t \in \mathcal{M}. \hspace{2.6cm}
\end{eqnarray}
\revI{The dynamics are expressed here in the embedding space $\mb{X}_t \in
\mathbb{R}^3$ and constrained to be within the curved surface
$\mathcal{M}$.}  We could also express these dynamics using local coordinate
charts on the surface as discussed in Section~\ref{sec:dynamics_ab}.  The $U$
is the potential energy, $\gamma$ the friction coefficient, $D =
{K_BT}/{\gamma}$ the diffusivity, and $K_B{T}$ the thermal energy with $K_B$
the Boltzmann constant and $T$ the temperature~\cite{Reichl1998}.  In terms of
the drift-diffusion tensors, $\mb{a}(\mb{X}_t) = {\gamma}^{-1} \gradX
U(\mb{X}), \;\; \mb{b}(\mb{X}_t) = \sqrt{2D} \mb{I}$.  

\subsection{Role of Drift in Mean First Passage Times: Double-Well Potential}
\begin{figure}[H]
\centering
\includegraphics[width=.99\textwidth]{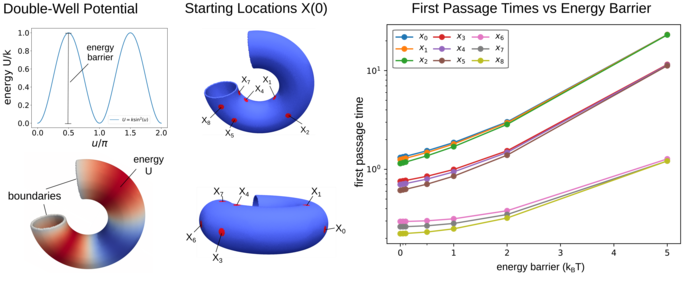}
  \caption{\revI{\textbf{FPTs with Drift-Diffusion Dynamics in Double-Well Potentials.}
\textit{(left)} Double-well potential $U$ on the surface of
equation~\ref{equ:double_well} influencing the drift of the stochastic
dynamics in~\ref{equ:langevin_surf}.  \textit{(middle}) The starting
locations for first passage times with $\mb{X}(0) = \mb{X}_i$.
\textit{(right)} First passage times from the selected starting points as the
strength of the energy barrier $k$ of the double-well potential is varied.
The locations $X_0,X_1,X_2$ start at the minimum of the potential energy and
exhibit the longest first passage times.  In this case the diffusion 
must surmount at least one of the energy barriers to reach the 
boundary.  As the energy barrier is
increased the first passage times increase, most rapidly when $k \gtrsim k_B
T$, and with an exponential trend in agreement with Kramer's
  theory~\cite{Gardiner1985}.}}
\label{fig:Drift_Study_Points_Chosen}
\label{fig:Drift_Study_First_Passage_Times}	
\end{figure}

Using our methods to capture drift in the surface dynamics, we study the role
that can be played by surface potentials in influencing stochastic trajectories
and the first passage time.  We consider the case of a double-well potential
$U$ on a surface which from the conservative force $F = -\nabla_X U$ in
equation~\ref{equ:langevin_surf} introduces a drift into the dynamics
influencing the mean first passage time.  We consider geometries $\Omega$
described by the truncated torus, parameterized by
\begin{eqnarray}
\mb{X}(u,v) = \begin{bmatrix} \cos(u)(.7 + .3\cos(v)), \;
 \sin(u)(.7 + .3\cos(v)), \;
.3\sin(v) \end{bmatrix}^T \nonumber, \,\,
u \in \left[ \frac{\pi}{4}, \frac{7\pi}{4}\right], \;\; v \in \left[0, 2 \pi \right).
\end{eqnarray} 
The double-well potential $U$ is generated using the parameterization $(u,v)$ as
\begin{equation}
\label{equ:double_well}
U(\mb{x}) = k \sin^2(u);\;\; \mb{x} = \mb{X}(u,v).
\end{equation}
We study how the mean first passage time changes as we vary the energy barrier $k = \tilde{k}\cdot k_BT $ with
$\tilde{k} \in 
\{0.0, 0.05, 0.1, 0.5, 1.0, 2.0, 5.0 \}.$
The case $\tilde{k} = 0.0$ serves as our baseline case with no drift.  The energy minimum occurs at $u = \pi$ and has energy barriers at $u = {\pi}/{2}$ and $u={3\pi}/{2}$ near the boundaries at $u = {\pi}/{4}$ and $u={7\pi}/{4}$.  We report the results of the first passage times starting at the specific locations $X_i$ given by
\begin{eqnarray}
\mb{X}_0: \left( u,v \right) = \left( \pi,0 \right) \nonumber \;
\mb{X}_1: \left( u,v \right) = \left( \pi , \frac{\pi}{2} \right) \nonumber \;
\mb{X}_2: \left( u,v \right) = \left( \pi , \pi \right) \nonumber \\
\mb{X}_3: \left( u,v \right) = \left( \frac{\pi}{2} , 0 \right) \nonumber \;
\mb{X}_4: \left( u,v \right) = \left( \frac{\pi}{2} , \frac{\pi}{2} \right) \nonumber \;
\mb{X}_5: \left( u,v \right) = \left( \frac{\pi}{2} , \pi \right) \nonumber \\
\mb{X}_6: \left( u,v \right) = \left( \frac{\pi}{3} , 0 \right) \nonumber \;
\mb{X}_7: \left( u,v \right) = \left( \frac{\pi}{3} , \frac{\pi}{2} \right) \nonumber \;
\mb{X}_8: \left( u,v \right) = \left( \frac{\pi}{3} , \pi \right). \nonumber 
\end{eqnarray}
We report our results in Figure~\ref{fig:Drift_Study_First_Passage_Times}.  

In our studies the starting locations $X_0,X_1,X_2$ are at the minimum of the
potential energy $U$.  From these locations each trajectory must surmount at
least one of the energy barriers to reach the boundary.  This results in the
largest first passage times.  We find as the energy barrier is increased $k
\gtrsim k_B T$ there is approximately exponential increase in the first passage
times.  This is in agreement with theory for first passage times involving
energy barriers based on asymptotic analysis using Kramer's
approximation~\cite{Gardiner1985}.  Our numerical methods allow for computing
the full solution $u(\mb{x})$ of equation~\ref{equ:u_x_stat}, including first
passage times, over a wide range of regimes $k \ll k_B T$, $k \sim k_B T$, and
$k \gg k_B T$, and when starting at locations that are not only at the
potential energy minimum, see Figure~\ref{fig:Drift_Study_Points_Chosen}.

\subsection{Role of Diffusion in Mean First Passage Times: Spatially-Dependent Diffusivities}

\begin{figure}[H]
\centering
\includegraphics[width=.99\textwidth]{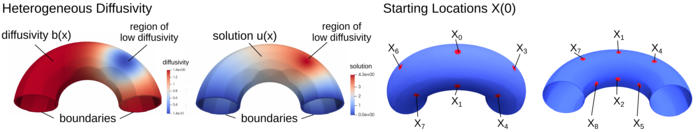} 
  \caption{\revI{\textbf{Spatially Dependent Diffusivities.}  \textit{(left)} The heterogeneous
spatially dependent diffusivity $b(\mb{x})$ in equation~\ref{equ:b_x}.
\textit{(left-middle)} The solution $u(\mb{x})$ for the first passage time
problem of equation~\ref{equ_surf_PDE_bk}.  \textit{(right-middle)},
\textit{(right)}  The starting locations $\mb{X}(0) = \mb{X}_i$.  We report
  results of studies in Figure~\ref{fig:Diff_Study_Mag_First_Passage_Times}.}  }
\label{fig:Diff_Study_Points_Chosen}
\end{figure}

\begin{figure}[H]
\centerline{\includegraphics[width=0.99\columnwidth]{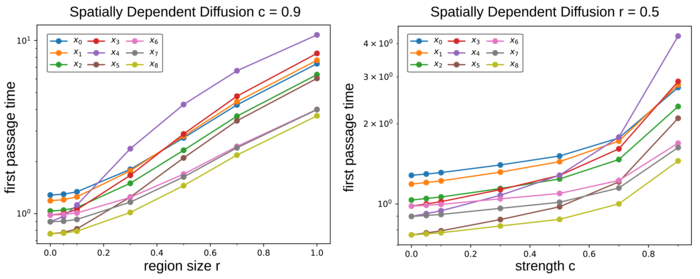}}
  \caption{\revI{\textbf{FPTs for Spatially Dependent Diffusivities.}
  \textit{(left)} The first passage times for the spatially dependent
  diffusivity $b(\mb{x})$ in equation~\ref{equ:b_x} when varying the size $r$
  of the region of low diffusivity.  \textit{(right)} The first passage times
  when increasing the strength $c$ of the diffusivity reduction in
  equation~\ref{equ:b_x}.}}
\label{fig:Diff_Study_Mag_First_Passage_Times}	
\label{fig:Diff_Stufy_Len_First_Passage_Times}	
\end{figure}

We study the role of spatially-dependent heterogeneous diffusivities $b = b(\mb{X})$ in first passage times.  We consider the surface geometry $\Omega$ obtained from the truncated torus
\begin{eqnarray}
\mb{X}(u,v) = \begin{bmatrix} 0.4\sin(v), \;
 \sin(u)(1 + 0.4\cos(v)) \nonumber, \;
-\cos(u)(1 + 0.4\cos(v))  \end{bmatrix}^T,\;\;
u \in \mathcal{U}(z = 0,v), v \in \left[0, 2 \pi \right). \nonumber
\end{eqnarray}
This geometry is obtained by slicing the torus at $z = 0$.  This defines the boundary $\partial \Omega$ and for the parameterization $(u,v)$ a range $\mathcal{U}$ for $u$ that depends on each value of $v$.  We consider 
the spatially-dependent diffusivity
\begin{eqnarray}
\label{equ:b_x}
\mb{b}(\mb{x}) = \sqrt{2D(\mb{x})}\mb{I} = \sqrt{2D_0} \left( 1 - c\, \exp\left({ \frac{ - \| \mb{x} - \mb{x}_c \|^2 }{r^2} }\right) \right) \mb{I}, \\ 
\mb{x}_c: \left( u_c, v_c \right) = \left( \frac{3 \pi}{4}, \frac{\pi}{2} \right), \;\;\;\; r \ge 0,\; 0 \le c < 1. 
\end{eqnarray}
This corresponds in equation~\ref{equ:langevin_surf} to 
$D = D(\mb{x}), \gamma = \gamma(\mb{x}),$ and zero drift 
with $U = 0$ so $\mb{a}(\mb{x}) = 0$.  The exponential is centered at $\mb{x}_c$ with decay over length-scale $r$ and influence amplitude $c$.  For $r$ sufficiently small the variations in diffusivity will occur primarily inside the domain $\Omega$ in a localized region $\|\mb{x} - \mb{x}_c\| \lesssim 3 r$.  The diffusivity elsewhere will be approximately constant $\sim D$.  
We consider for our initial starting locations
\begin{eqnarray}
\mb{X}_0: \left( u,v \right) = \left( \pi,0 \right) \nonumber \;
\mb{X}_1: \left( u,v \right) = \left( \pi , \frac{\pi}{2} \right) \nonumber \;
\mb{X}_2: \left( u,v \right) = \left( \pi , \pi \right) \nonumber \\
\mb{X}_3: \left( u,v \right) = \left( \frac{3\pi}{4} , 0 \right) \nonumber \;
\mb{X}_4: \left( u,v \right) = \left( \frac{3\pi}{4} , \frac{\pi}{2} \right) \nonumber \;
\mb{X}_5: \left( u,v \right) = \left( \frac{3\pi}{4} , \pi \right) \nonumber \\
\mb{X}_6: \left( u,v \right) = \left( \frac{5\pi}{4} , 0 \right) \nonumber \;
\mb{X}_7: \left( u,v \right) = \left( \frac{5\pi}{4} , \frac{\pi}{2} \right) \nonumber \;
\mb{X}_8: \left( u,v \right) = \left( \frac{5\pi}{4} , \pi \right). \nonumber
\end{eqnarray}
We study the role of spatially-dependent diffusivity in two cases: (i) when the
depth $c$ is varied while leaving the area scale $r$ fixed, and (ii) when the
area scale $r$ is varied while leaving the depth $c$ fixed.  We use  throughout
the baseline parameters $c = 0.9$ and $r = 0.5$.  We show the geometry and
solution for a typical spatially dependent diffusivity in
Figure~\ref{fig:Diff_Study_Points_Chosen}.  We report our results in
Figure~\ref{fig:Diff_Study_Mag_First_Passage_Times}.

We find varying the spatial extent $r$ of the diffusivity plays a particularly
strong role in the first passage times.  The underlying mechanisms are
different than the double-well potential.  When a trajectory approaches
$\mb{x}_c$ the motion of the stochastic process slows down significantly as a
result of the smaller diffusivity.  There is no barrier for $\mb{X}_t$ to
approach such regions but once within this region exhibits a type of temporary
trapping behavior from the slow diffusion.  We see from
Figure~\ref{fig:Diff_Study_Points_Chosen}, such a region influences the first
passage time over a much larger range than the direct variations in $b(\mb{x})$
given the high probability of a large fraction of trajectories encountering
this trapping region even when a modest distance away.  We see increases in $r$
have a strong influence on increasing the first passage times, see
Figure~\ref{fig:Diff_Study_Mag_First_Passage_Times}.  We also see in the limit
that the localized diffusivity $b(\mb{x})$ approaches zero, even a relatively
small probability of encountering this region can result in a large first
passage time.  We see this increasing influence as $c$ approaches one, see
Figure~\ref{fig:Diff_Study_Mag_First_Passage_Times}.

\subsection{Role of Geometry in Mean First Passage Times: Neck-Shaped Domains}
\begin{figure}[H]
\centering
\includegraphics[width=.99\textwidth]{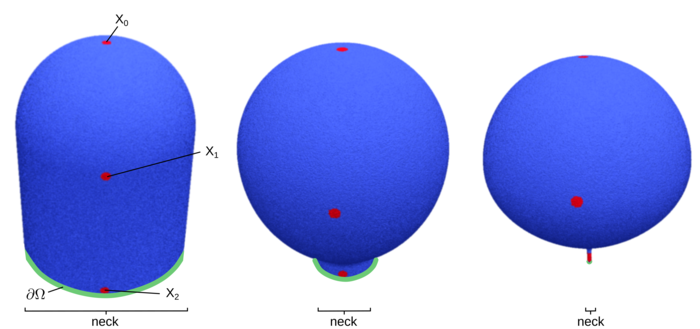} 
\hfill 
  \caption{\revI{\textbf{Neck Geometries for First Passage Time Studies.}
  For the FPT studies when varying the geometry of the neck region, we
show the select starting locations $\mb{X}(0) = \mb{X}_i$ \textit{(red
points)}.  On the left we show the largest neck domain which does not present
a geometric bottleneck for reaching the $\partial \Omega$ at $z = 0$.  
  On the right, we show the
smallest neck region which presents a significant geometric bottleneck
potentially inhibiting diffusion in reaching $z = 0$.  The geometries used in
these studies were designed so that as the neck width is varied the geodesic
distance from the starting points remains constant to reach the boundary $\partial
  \Omega$ at
$z =0$.}} \label{fig:Neck_Study_Points_Chosen}
\end{figure}

\begin{figure}[H]
\centerline{\includegraphics[width=0.99\columnwidth]{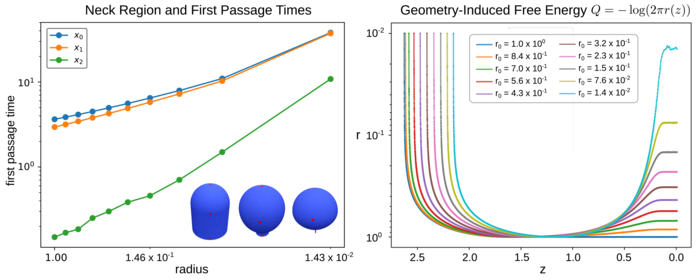}}
  \caption{\revI{\textbf{FPTs Dependence on the Neck Geometry.}
  \textit{(left)} The first passage time as the neck radius $r$ becomes
  small.  We see as the neck becomes narrower (towards the right) the geometry
  of the manifold restricts the diffusion and reaching the boundary becomes
  more difficult.  The $r$ indicates the narrowest radius of the shape.
  \textit{(right)} The shape acts to give a geometrically induced free energy
  $Q(z) = -log(2\pi r(z))$ for the dynamics when projecting to the $z$-axis.
  We see as the neck radius $r$ narrows (towards the right) the geometry
  creates an increasingly large free-energy barrier and increases the
  first-passage times. } }
\label{fig:Neck_Study_Energy_Barrier}
\label{fig:Neck_Study_First_Passage_Times}
\end{figure}

\revI{We show how our methods can be used to investigate how the shape of the
curved surface influences first passage times.}  In these studies, we
consider purely diffusive dynamics.  This corresponds to the case with $U = 0$
with zero drift $\mb{a}(\mb{x}) = 0$ and a spatially homogeneous diffusivity
$\mb{b}(\mb{x}) = \sqrt{2D_0} I$ in equation~\ref{equ:langevin_surf}.  For the
geometry we use a surface of revolution with an adjustable shape that forms a
neck region near the bottom at $z = 0$, see
Figure~\ref{fig:Neck_Study_Points_Chosen}.  This geometry is generated by
considering first a cylinder of radius $r_0$ and height $h = 0.05$ which is
capped at the top by a sphere.  We then use a radial profile $r(z)$ to connect
the spherical cap with the bottom of the cylinder of radius $r_0$ at $z = 0$.
For this purpose, we use a bump function $r(z) = (1- 0.4b) + 0.4b \exp\left({1
- \frac{b^2}{b^2 - b + .05 -z)^2}}\right)$.  For $z \in [0.05, b + 0.05]$, we
choose $b$ so that this has arc-length ${\pi}/{2}$, and the final radius $r_0$
at $z = 0$. This serves to smoothly connect the unit hemisphere cap to the
bottom of the cylinder with radius $r_0$.  This ensures the geodesic distance
from $X(0)=X_i$ to the boundary is always the same as we vary the shapes, see
Figure~\ref{fig:Neck_Study_Points_Chosen}.  We report our results in
Figure~\ref{fig:Neck_Study_First_Passage_Times}.

We find as the neck region becomes smaller it acts as a hindrance for
trajectories to reach the boundary and first passage times increase.  From the
log-log plot we see that the first passage time appears over many radii to
follow a power-law trend $\bar{\tau} \sim r^{\alpha}$ with values for $\alpha$
respectively for $X(0) = X_0,X_1,X_2$ having $\alpha \approx 0.427,
0.484,0.893$.  As the radius tends to $0$ we see the first passage time
diverges.  We see for points further from the boundary the first passage times
are longer as would be expected.  However, as the neck region becomes small, we
see smaller differences occur between points $X_0,X_1$ starting at points above
the neck region, compared to the point $X_3$ starting within the neck region.
This indicates that as the radius shrinks the neck region increasingly acts as
a hindrance for reaching the boundary.  Interesting, we also see that starting
at $X_3$ also has FPTs that significantly increase, since as $r_0 \rightarrow
0$ an increasingly large fraction of trajectories will leave the neck region
before encountering the boundary and then must also overcome the hindrance
similar to starting at $X_0,X_1$.  We can characterize this mechanism by using
a reaction-coordinate and notion of free energy (entropy contributions) arising
from the constricting geometry.

Given the radial symmetry we can project the stochastic dynamics onto the
$z$-component as an effective reaction-coordinate for reaching the boundary at
$z = 0$.  This provides a statistical mechanics model~\cite{Reichl1998} where
the shape of the manifold gives an effective free energy having entropic
contributions proportional to $Q(z) = -\log(2\pi r(z))$.  This contributes to
an effective mean-force for $Z_t$ yielding the drift $a(z) = -\partial
Q/\partial z = -Q'(z)$.  This suggests the radially-averaged stochastic
dynamics $ dZ_t = -Q'(Z_t)dt + C\,dW_t.  $ From this perspective, the
drift-diffusion process $Z_t$ must cross the geometrically-induced energy
barrier $Q(z)$ to reach the boundary at $z = 0$.  Similar to our results on
double-well potentials, this could be studied analytically with asymptotic
Kramer's theory~\cite{Gardiner1985}.  Our numerical methods allow for capturing
directly these effects and we see as the neck narrows this significantly
increases the first passage times.  This correlates well with how the effective
free energy barrier for $Z_t$ increases with the shape change, see
Figure~\ref{fig:Neck_Study_Energy_Barrier}.  Our results show how geometric
effects captured by our methods can provide mechanisms that influence
significantly the behaviors of stochastic processes and their first passage
times.
 
\section*{Conclusions} 
\revI{We have developed numerical methods for computing
the mean First Passage Times and related path-dependent 
statistics of 
stochastic processes on surfaces of general shape and topology.
We formulated descriptions of Brownian motion and general
drift-diffusion processes on surfaces.
Using Dynkin's formula and Backward-Kolmogorov equations
of these processes, we formulated and 
solved associated elliptic PDE boundary value problems 
on curved surfaces.  We developed our numerical
methods using Generalized Moving Least Squares (GMLS)
to approximate the local surface geometry and 
the action of surface differential operators. 
Using this discretization approach we introduced 
collocation methods and solvers for the associated 
linear systems of equations.  For a variety of 
surface shapes, we showed that our methods converge 
with high-order accuracy in capturing the geometry 
and PDE solutions.  For mean First Passage Times, 
we showed how our methods can be used to
investigate the roles of the surface geometry, drift
dynamics, and spatially dependent diffusivities.  The solvers 
and approaches we have developed also can be used more generally
to compute other path-dependent statistics 
and solutions to elliptic boundary-value 
problems on surfaces of general shape.}

\section*{Acknowledgements}
The authors would like to acknowledge support from research grants DOE ASCR
PhILMS DE-SC0019246 and NSF DMS-1616353.


\printbibliography


\appendix
\addcontentsline{toc}{section}{Appendices}
\section*{Appendix}

\section{Monge-Gauge Surface Parameterization}
\label{appendix:monge_gauge_param}
In the Monge-Gauge we parameterize locally a smooth surface in terms of the
tangent plane coordinates $u,v$ and the height of the surface above this point
as the function $h(u,v)$.  This gives the embedding map
\begin{eqnarray}
\label{equ:manifoldParam}
\mb{x}(u,v) = \bs{\sigma}(u, v) =  (u, v, \h(u,v)).
\end{eqnarray}
We can use the Monge-Gauge equation~\ref{equ:manifoldParam} to derive explicit
expressions for geometric quantities.  The derivatives of $\bsy{\sigma}$
provide a basis $\partial_u, \partial_v$ for the tangent space as
\begin{eqnarray}
\label{eqn:sigmaU}
\partial_u = \bs{\sigma}_{u}(u, v) & = & (1, 0, \h_u(u, v) ) \\
\label{eqn:sigmaV}
\partial_v = \bs{\sigma}_{v}(u, v) & = & (0, 1, \h_v(u, v) ).
\end{eqnarray}
The first fundamental form $\mathbf{I}$ (metric tensor) and its inverse $\mathbf{I}^{-1}$ (inverse tensor) are given by 
\begin{align}
\mathbf{I} = \begin{bmatrix}
g_{uu} & g_{uv} \\
g_{vu} & g_{vv}
\end{bmatrix} = \begin{bmatrix}
\bs{\sigma}_{u} \cdot \bs{\sigma}_{u} & \bs{\sigma}_{u} \cdot \bs{\sigma}_{v} \\
\bs{\sigma}_{v} \cdot \bs{\sigma}_{u} & \bs{\sigma}_{v} \cdot \bs{\sigma}_{v}
\end{bmatrix} =  \begin{bmatrix}
1 + \h_u(u,v)^2 & \h_u \h_v(u,v) \\
\h_u(u,v) \h_v(u,v) & 1 + \h_v(u,v)^2
\end{bmatrix}. \label{equ:I_mat}
\end{align}
and 
\begin{align}
\mathbf{I}^{-1} = \begin{bmatrix}
g^{uu} & g^{uv} \\
g^{vu} & g^{vv}
\end{bmatrix} = \frac{1}{ 1 + \h_u^2 + \h_v^2 }\begin{bmatrix}
1+ \h_{v}^2 & -\h_{u}\h_{v} \\
-\h_{u}\h_{v} & 1+ \h_{u}^2
\end{bmatrix} . \label{equ:I_inv_mat} 
\end{align}

We use throughout the notation for the metric tensor $\mb{g} = \mathbf{I}$  interchangeably.  For notational convenience, we use the tensor notation for the metric tensor $g_{ij}$ and for its inverse $g^{ij}$.  These correspond to the first fundamental form and its inverse as
\begin{eqnarray}
g_{ij} = \left[\mathbf{I}\right]_{i,j}, \hspace{0.3cm} g^{ij} = \left[ \mathbf{I}^{-1}\right]_{i,j}.
\end{eqnarray}
For the metric tensor $\mb{g}$, we also use the notation $|g| = \det(\mb{g})$ and 
have that
\begin{eqnarray}
\sqrt{|g|} = \sqrt{\det(\mathbf{I})} 
= \sqrt{1 + \h_u^2 + \h_v^2} = \| \mb{\sigma}_{u}(u, v) \times \mb{\sigma}_{v}(u, v) \|.
\label{met_fac_def}
\end{eqnarray}
The provides the local area element as $dA_{u,v} = \sqrt{|g|}du dv$. 

\revI{The Christoffel symbols of the second kind were used in
our derivations of Brownian motion and more general stochastic
drift-diffusion dynamics on surfaces. 
In the general case, the Christoffel symbols are given by
\begin{eqnarray}
  \Gamma^k_{ij} = \partial_j \bs{\sigma}_i \cdot \bs{\sigma}_k.
\end{eqnarray}
In the case of the Monge-Gauge $\bs{\sigma}(u,v) = (u,v,\h(u,v)$, this can be expressed as
\begin{equation}
\Gamma^k_{ij} = \frac{ \h_{ij} \h_k }{ \sqrt{|g|} } = \frac{ \h_{ij} \h_k }{  \sqrt{1 + \h_u^2 + \h_v^2} }.
\end{equation}
To compute quantities associated with 
curvature of the manifold, we construct 
the Weingarten map~\cite{Pressley2001}.  This can be expressed as
\begin{eqnarray}
\mb{W} = \mb{I}^{-1} \mb{II}.
\end{eqnarray}
The Gaussian curvature $K$ can be expressed in the Monge-Gauge as
\begin{eqnarray}
K(u,v) = \det\left(\mb{W}(u,v)\right) = \frac{\h_{uu}\h_{vv} 
- \h_{uv}^2}{(1 + \h_u^2 + \h_v^2)^2}.
\end{eqnarray}
For further discussions on the differential
geometry of manifolds and computation of surface operators 
see~\cite{Pressley2001,Abraham1988,SpivakDiffGeo1999,
GMLSStokes2020}.
}


\section{Sampling Resolution of the Manifolds}
\label{appendix:sampling_res_manifolds} We provide a summary of the sampling
resolution $h$ used for each of the manifolds in
Table~\ref{table:sampling_manifolds}.  We refer to $h$ as the \textit{target
fill distance}.  For each of the manifolds and target spacing $h$, we achieve a
nearly uniform collection of the points (quasi-uniform samplings) using
DistMesh~\cite{DistMesh2004}.  In practice,  we have found this yields a point
spacing with neighbor distances varying by about $\approx \pm 30 \%$ relative
to the target distance $h$.  We summarize for each of the manifolds how this
relates to the number of sample points $n$ in
Table~\ref{table:sampling_manifolds}.  Robustness studies for related solvers
have been carried out when applying perturbations to such point samples
in~\cite{GMLSStokes2020}.

\begin{center}
 \begin{tabular}{||c | c c |c c |c c |c c ||} 
 \hline
Refinement Level & \textbf{A}: h \hspace{0.22cm} & n & \textbf{B}: h \hspace{0.22cm} & n & \textbf{C}: h \hspace{0.22cm} & n & \textbf{D}: h \hspace{0.22cm} & n \\ 
 \hline\hline
 $1$ & .1 & 2350 & .1 & 2306 & .1 & 2002 & .08 & 1912 \\ 
 \hline
 $2$ & .05 & 9566 & .05 & 9206 & .05 & 7998 & .04 & 7478 \\
 \hline
 $3$ & .025 & 38486 & .025 & 36854 & .025 & 31898 & .02 & 29494 \\
 \hline
 $4$ & .0125 & 154182 & .0125 & 147634 & .0125 & 127346 & .01 & 118942 \\
 \hline
\end{tabular}
\captionof{table}{Sampling Resolution for each of the Manifolds A--D.  Relation
between the target distance $h$ and the number of sample points $n$ used for
each of the manifolds.  In each case, the neighbor distances between the
points sampled were within $\approx \pm 30\%$ of the target distance $h$.  } 
\label{table:sampling_manifolds}
\end{center}

\end{document}